\documentclass[draft, a4paper]{article}
\usepackage{amsmath, amsfonts, amsthm, amssymb}
\numberwithin{equation}{section}

\newcommand{\Q}{\mathbb{Q}}

\newcommand{\R}{\mathbb{R}}
\newcommand{\C}{\mathbb{C}}
\newcommand{\N}{\mathbb{N}}
\newcommand{\Z}{\mathbb{Z}}

\newtheorem{thm}{Theorem}
\newtheorem{lem}{Lemma}

\renewcommand{\mod}[1]{\hspace{-2.9mm}\pmod{#1}}
\newcommand{\x}{{\bf x}}

\newcommand{\rom}{\mathrm}
\newcommand{\bfP}{\mathbb{P}}
\newcommand{\A}{\mathbb{A}}
\newcommand{\ov}{\overline}
\newcommand{\ma}{\mathbf}
\newcommand{\ben}{\begin{enumerate}}
\newcommand{\een}{\end{enumerate}}
\newcommand{\eit}{\begin{itemize}}

\newcommand{\ve}{\varepsilon}

\newcommand{\mcal}{\mathcal}
\newcommand{\lab}{\label}
\newcommand{\al}{\alpha}

\newcommand{\colt}[2]{\genfrac{}{}{0pt}{1}{#1}{#2}}

\renewcommand{\leq}{\leqslant}
\renewcommand{\geq}{\geqslant}

\newcommand\HH{\mathcal{H}}

\newcommand\NU{{N_{U,H}}}
\newcommand{\Esix}{{\mathbf E}_6}
\newcommand{\Dfour}{{\mathbf D}_4}
\newcommand{\Dfive}{{\mathbf D}_5}
\newcommand{\tS}{{\widetilde S}}
\newcommand{\tV}{{\widetilde V}}

\newcommand{\nub}{N_{U,H}(B)}

\DeclareMathOperator{\Pic}{Pic}

\DeclareMathOperator{\Proj}{Proj}
\DeclareMathOperator{\hcf}{gcd}

\newtheorem*{con1}{Conjecture A}
\newtheorem*{con2}{Conjecture B}
\newtheorem*{con3}{Conjecture C}

\theoremstyle{definition}
\newtheorem*{ack}{Acknowledgements}

\begin{document}

\title{An overview of Manin's conjecture \\for del Pezzo surfaces}
\author{T.D. Browning\\
\small{\emph{School of Mathematics, 
Bristol University, Bristol BS8 1TW}}\\
\small{t.d.browning@bristol.ac.uk}}

\date{}
\maketitle

\section{Introduction}\lab{intro}

A fundamental theme in mathematics is the study of integer solutions
to Diophantine equations, or equivalently, the study of rational 
points on projective algebraic varieties.  Let $V \subset \bfP^n$ be a
projective variety that is cut out by a finite system of
homogeneous equations defined over $\Q$.  
Then there are a number of basic questions that can be asked about the
set $V(\Q):=V\cap \bfP^n(\Q)$ of rational points on $V$:
when is $V(\Q)$ non-empty? how large is $V(\Q)$ when it is non-empty?
This paper aims to survey the second question, for a large class of
varieties $V$ for which one expects $V(\Q)$ to be Zariski dense in $V$.

To make sense of this it is convenient to define the {\em height} of
a projective rational point $x=[x_0,\ldots, x_{n}] \in \bfP^{n}(\Q)$ to be
$H(x):=\|\x\|$, for any norm $\|\cdot\|$ on $\R^{n+1}$,  
provided that $\x=(x_0, \ldots, x_{n}) \in \Z^{n+1}$ and
$\hcf(x_0,\ldots,x_{n})=1$. Throughout this work we shall work with
the height metrized by the choice of norm $|\x|:=\max_{0\leq i \leq
  n}|x_i|$.  Given a suitable Zariski open subset $U\subseteq V$, the
goal is then to study the quantity
\begin{equation}\lab{count}
N_{U,H}(B):= \#\{ x \in U(\Q):  \,H(x) \leq B\}, 
\end{equation}
as $B \rightarrow \infty$. It is natural to question whether the 
asymptotic behaviour of $N_{U,H}(B)$ can be related to the geometry of
$V$, for suitable open subsets $U\subseteq V$.
Around 1989 Manin initiated a program to do exactly this for
varieties with ample anticanonical divisor \cite{f-m-t}.  
Suppose for simplicity that $V\subset \bfP^n$ is a non-singular complete
intersection, with $V=W_1\cap \cdots \cap W_t$ for hypersurfaces $W_i \subset \bfP^n$ of degree $d_i$.
Since $V$ is assumed to be Fano, its Picard group is a finitely
generated free $\Z$-module, and we denote its rank by $\rho_V$.  
Then in this setting the Manin conjecture takes the following shape
\cite[Conjecture $C'$]{b-m}.

\begin{con1}
Suppose that $d_1+\cdots+d_t\leq n$.  Then there exists a Zariski open
subset $U \subseteq V$ and a non-negative constant $c_{V, H}$ such that 
\begin{equation}\lab{c1}
N_{U,H}(B) = c_{V,H}B^{n+1-d_1-\cdots-d_t} (\log B)^{\rho_V-1}\big(1+o(1)\big), 
\end{equation}
as $B\rightarrow \infty$.
\end{con1}

It should be noted that there are simple heuristic arguments that
support the value of the exponents of $B$ and $\log B$ appearing in
the conjecture.  The constant $c_{V,H}$ has also received
a conjectural interpretation  at the hands of Peyre \cite{MR1340296},
and this has been generalised to cover certain other cases by Batyrev and
Tschinkel~\cite{b-t}, and  Salberger~\cite{MR1679841}.  
In fact whenever we refer to the Manin conjecture we shall henceforth
mean that the value of the constant $c_{V,H}$ should agree with the
prediction of Peyre et al. With this in mind, the 
Manin conjecture can be extended to cover certain
other Fano varieties $V$ which are not necessarily complete intersections, nor 
non-singular. For the former one simply takes the exponent of $B$ to be the infimum
of all $a/b \in \Q$ such that $b>0$ and $aH+bK_V$ is linearly
equivalent to an effective divisor, where $K_V \in \rom{Div} (V)$ is a
canonical divisor and $H \in \rom{Div} (V)$ is a hyperplane section.
For the latter, if $V$ has only rational double points, then one may
apply the conjecture to a minimal desingularisation $\tV$ of
$V$, and then employ the functoriality of heights.  
A discussion of these more general versions of the conjecture
can be found in the survey of Tschinkel \cite{t1}.
The purpose of this note is to give an overview of our progress in the
case that $V$ is a suitable Fano variety of dimension $2$.

A non-singular surface $S \subset \bfP^d$ of degree $d$, with very ample anticanonical divisor
$-K_S$, is known as a {\em del Pezzo surface of degree $d$}.  Their geometry has been expounded by Manin
\cite{manin-book}, for example.  It is well-known that such
surfaces $S$ arise either as the quadratic Veronese embedding of a quadric
in $\bfP^3$, which is a del Pezzo surface of degree $8$ in
$\bfP^8$ (isomorphic to $\bfP^1\times \bfP^1$), or as the blow-up of $\bfP^2$ along $9-d$ points in general
position, in which case the degree of $S$ satisfies $3
\leq d \leq 9$.  In terms of the expected asymptotic formula for
$N_{U,H}(B)$ for a suitable open subset $U \subseteq S$, the exponent
of $B$ is $1$, and the exponent of
$\log B$ is at most $9-d$, since the geometric Picard group $\Pic(S
\otimes_\Q \overline{\Q})$ has rank $10-d$.
An old result of Segre ensures that the set $S(\Q)$ of rational points
on $S$ is Zariski dense as soon as it is non-empty.
Moreover, when $3 \leq d\leq 8$ there are certain so-called
{\em accumulating subvarieties} contained in $S$ which may dominate
the behaviour of the counting function $N_{S,H}(B)$.  These are the
possible lines contained in $S$, that correspond to the exceptional
divisors arising from the process of blowing up the projective 
plane along the relevant collection of points.  Now it is an easy
exercise to check that 
$$
N_{\bfP^1,H}(B)=\frac{12}{\pi^2}B^2\big(1+o(1)\big),
$$
as $B \rightarrow \infty$, so that $N_{V,H}(B)\gg_{V} B^2$ for any
geometrically integral surface $V\subset \bfP^n$ that contains a line
defined over $\Q$.  However, if $U\subseteq V$ is defined to be the
Zariski open subset formed by deleting all of the lines from $V$ then
it follows from combining an estimate
of Heath-Brown \cite[Theorem $6$]{annal} with a birational projection
argument due to Salberger \cite[\S $8$]{s}, that 
$N_{U,H}(B)=o_V(B^{2})$.

Returning to the setting of del Pezzo surfaces $S \subset \bfP^d$ of
degree $d$, it turns out that there are no
exceptional divisors when  $d=9$, or when $d=8$ and $S$ is isomorphic
to $\bfP^1\times \bfP^1$, in which case we study
$N_{S,H}(B)$. When $3 \leq d\leq 7$, or when $d=8$ and $S$ is not
isomorphic to $\bfP^1\times \bfP^1$, there are a finite number of such
divisors, each producing a line in $S$. In these cases we study $N_{U,H}(B)$ for the
open subset $U$ formed by deleting all of the lines from $S$.
We now proceed to review the progress that has been made towards the Manin
conjecture for del Pezzo surfaces.  
In doing so we have split our discussion according to the degree of the surface,
and it will become apparent that the problem of estimating
$N_{U,H}(B)$ becomes harder as the degree decreases.

\subsection{Del Pezzo surfaces of degree $\geq 5$}

It turns out that the del Pezzo surfaces $S$ of degree $d\geq 6$ are toric,
in the sense that they contain the torus $\mathbb{G}_m^{2}$ of algebraic
groups as a dense open subset, whose natural action on itself extends to all of $S$.  
Thus the Manin conjecture for such surfaces is a special case of the
more general work due to Batyrev and Tschinkel \cite{MR1620682}, that
establishes this conjecture for all toric varieties.  One may compare this result with the work of 
la Bret\`eche \cite{dlB-toric} and Salberger \cite{MR1679841}, who
both establish the conjecture for toric varieties defined over $\Q$,
and also the work of Peyre \cite{MR1340296}, who handles a number
of special cases.

For non-singular del Pezzo surfaces $S \subset \bfP^5$ of degree~$5$,
the situation is rather less satisfactory.  In fact there are very few
instances for which the Manin conjecture has been
established.  The most significant of these is due to la Bret\`eche
\cite{MR2003m:14033}, who has proved the conjecture when all of the $10$
exceptional divisors are defined over $\Q$.  In this case we say that
the surface is {\em split}.   Let $S_0$ be the surface 
obtained by blowing up $\bfP^2$ along the four points
$$
p_1=[1,0,0], \quad p_2=[0,1,0], \quad p_3=[0,0,1], \quad p_4=[1,1,1],
$$  
and let $U_0\subset S_0$ denote the corresponding open subset formed
by deleting the lines from $S_0$.  Then $\Pic (S_0)$ has rank $5$,
since $S_0$ is split, and la Bret\`eche obtains the following result.

\begin{thm}\lab{5-1}
Let $B \geq 3$.  Then  there exists a constant $c_0>0$ such that
$$
N_{U_0,H}(B) = c_{0}B (\log B)^{4}\Big(1+O\Big(\frac{1}{\log\log B}\Big)\Big).
$$
\end{thm}

We shall return to the proof of this result below.  The other 
major achievement in the setting of quintic del Pezzo surfaces is a
result of la Bret\`eche and Fouvry \cite{MR2099200}.  Here the
Manin conjecture is established for the surface obtained by
blowing up $\bfP^2$ along four points in general position, two of
which are defined over $\Q$ and two of which are conjugate over $\Q(i)$.
In related work, Browning  \cite{dp5} has obtained upper bounds for 
$N_{U,H}(B)$ that agree with the Manin predication for several other
del Pezzo surfaces of degree $5$.

\subsection{Del Pezzo surfaces of degree $4$}

A quartic del Pezzo surface $S \subset \bfP^4$, that is defined over $\Q$, can be
recognised as the zero locus of a suitable pair of quadratic forms 
$Q_1,Q_2 \in \Z[x_0,\ldots,x_4]$.  Then 
$S=\Proj(\Q[x_0,\ldots,x_4]/(Q_1,Q_2))$ 
is the complete intersection of the hypersurfaces $Q_1=0$ and $Q_2=0$ in
$\bfP^4$. When $S$ is non-singular \eqref{c1} predicts 
the existence of a constant $c_{S,H} \geq 0$ such that
\begin{equation}\lab{c2}
\nub= c_{S,H} B (\log B)^{\rho_S-1}\big(1+o(1)\big),
\end{equation}
as $B \rightarrow \infty$, where $\rho_S=\rom{rk} \Pic(S) \leq 6$ and
$U\subset S$ is obtained by deleting the $16$ lines from $S$.
In this setting the best result available is due to
Salberger.  In work communicated 
at the conference {\em Higher dimensional varieties and rational
  points} at Budapest in 2001, he establishes the
estimate $\nub =O_{\ve,S}(B^{1+\ve})$ for any $\ve>0$, provided that 
the surface contains a conic defined over $\Q$.  In fact an examination of
Salberger's approach, which is based upon fibering the surface into a
family of conics, reveals that it would be straightforward to replace the factor
$B^\ve$ by $(\log B)^A$ for a large constant $A$.  
It would be interesting to find examples of surfaces $S$ for which the
exponent $A$ could be reduced to the expected quantity $\rho_S-1$.

It emerges that much more can be said if one permits $S$ to contain 
isolated singularities.  For the remainder of this section let 
$S \subset \bfP^4$ be a geometrically integral intersection of two
quadric hypersurfaces, which has only isolated singularities and is
not a cone. Then $S$ contains only rational double points
(see Wall \cite{wall}, for example), thereby ensuring that there exists a
unique minimal desingularisation $\pi: \tS\rightarrow S$ of the
surface, such that $K_{\tS}=\pi^* K_S$.
In particular it follows that the asymptotic formula \eqref{c2} is
still expected to hold, with $\rho_S$ now taken to be the rank of the Picard group
of $\tS$, and $U\subset S$ obtained  by deleting all of the lines from $S$.
The classification of such surfaces $S$ is rather
classical, and can be found in the work of Hodge and Pedoe \cite[Book IV, \S
XIII.11]{h-p}, for example.  The notation used there is rather
old-fashioned however, and makes it difficult to follow.  
Let $S=\Proj(\Q[\x]/(Q_1,Q_2))$ be as above. Then it turns out that up to isomorphism over $\overline{\Q}$, there are
$15$ possible singularity types for $S$, each categorised by the
{\em extended Dynkin diagram}. This is the Dynkin diagram that describes
the intersection behaviour of the exceptional divisors and the
transforms of the lines on the minimal desingularisation $\tS$ of $S$.
Of course, if one is interested in a classification over the ground
field $\Q$, then many more singularity types can occur (see Lipman
\cite{lipman}, for example).
Over $\overline{\Q}$, Coray and Tsfasman  \cite[Proposition 6.1]{c-t} have calculated the
extended Dynkin diagrams for all of the $15$ types, and this
information allows us to write down a list of surfaces that typify
each possibility, together with their singularity type and the number
of lines that they contain.   
The author is grateful to Ulrich Derenthal for helping to prepare the
following table, which lists examples of surfaces 
$S=\Proj(\Q[\x]/(Q_1,Q_2))$ that illustrate the possible types.

\begin{table}[!h]
\begin{center}
\begin{tabular}{|c|c|c|c|c|}
\hline
type & $Q_1(\x)$ & $Q_2(\x)$ & $\#$ lines & singularity \\
\hline
\hline
\texttt{i} & $x_0x_1+x_2x_3$ & $x_0x_3+x_1x_2+x_2x_4+x_3x_4$ & $12$ & $\textbf{A}_1$\\
\texttt{ii} & $x_0x_1+x_2x_3$ & $x_0x_3+x_1x_2+x_2x_4+x_4^2$ & $9$ & $2\textbf{A}_1$\\
\texttt{iii} & $x_0x_1+x_2^2$ & $x_0x_2+x_1x_2+x_3x_4$ & $8$ & $2\textbf{A}_1$\\
\texttt{iv} & $x_0x_1+x_2x_3$ &
$x_2x_3+x_4(x_0+x_1+x_2-x_3)$ & $8$ & $\textbf{A}_2$\\
\texttt{v} & $x_0x_1+x_2^2$ & $x_1x_2+x_2^2 +x_3x_4$   & $6$ & $3\textbf{A}_1$ \\
\texttt{vi} & $x_0x_1+x_2x_3$ & $x_1^2+x_2^2+x_3x_4$   & $6$ & $\textbf{A}_1+\textbf{A}_2$ \\
\texttt{vii} & $x_0x_1+x_2x_3$ &  $x_1x_3+x_2^2+x_4^2$ & $5$& $\textbf{A}_3$ \\
\texttt{viii} & $x_0x_1+x_2^2$ & $(x_0+x_1)^2+x_2x_4+x_3^2$   & $4$& $\textbf{A}_3$ \\
\texttt{ix} & $x_0x_1+x_2^2$ & $x_2^2+x_3x_4$   & $4$ & $4\textbf{A}_1$ \\
\texttt{x} & $x_0x_1+x_2^2$ & $x_1x_2+x_3x_4$   & $4$ & $2\textbf{A}_1+\textbf{A}_2$ \\
\texttt{xi} & $x_0x_1+x_2^2$ & $x_0^2+x_2x_4+x_3^2$   & $3$ & $\textbf{A}_1+\textbf{A}_3$ \\
\texttt{xii} & $x_0x_1+x_2x_3$ & $x_0x_4+x_1x_3+x_2^2$   & $3$ & $\textbf{A}_4$ \\
\texttt{xiii} & $x_0x_1+x_2^2$ & $x_0^2+x_1x_4+x_3^2$   & $2$ & $\textbf{D}_4$ \\
\texttt{xiv} & $x_0x_1+x_2^2$ & $x_0^2+x_3x_4$   & $2$ & $2\textbf{A}_1+\textbf{A}_3$ \\
\texttt{xv} & $x_0x_1+x_2^2$  &$x_0x_4+x_1x_2+x_3^2$   & $1$ & $\textbf{D}_5$ \\
\hline
\end{tabular}
\end{center}
\end{table}

Let $\tS$ denote the minimal desingularisation of any surface $S$ from
the table, and let
$\rho_S$ denote the rank of the Picard group of  $\tS$.  
Then it is natural to try and establish \eqref{c2} for such surfaces
$S$.  Several of the surfaces are actually special cases of
varieties for which the Manin conjecture is already known to hold.  
Thus we have seen above that it has been established for toric
varieties, and it can be checked that the 
surfaces representing types \texttt{ix}, \texttt{x}, \texttt{xiv} are
all equivariant compactifications of $\mathbb{G}_m^2$, and so are
toric.  Hence \eqref{c2} holds for
these particular surfaces.  Similarly it has been shown by 
Chambert-Loir and Tschinkel \cite{ct} that the Manin conjecture
is true for equivariant compactifications of vector groups.
Although identifying such surfaces in the table is not entirely routine,
it transpires that the $\Dfive$ surface representing  
type \texttt{xv} is an equivariant compactification of $\mathbb{G}_a^2$.
Per Salberger has raised the question of whether there exist singular del
Pezzo surfaces of degree $4$ that arise as equivariant
compactifications of $\mathbb{G}_a \times \mathbb{G}_m$,  
but  that are not already equivariant compactifications of
$\mathbb{G}_a^2$ or $\mathbb{G}_m^2$.
This is a
natural class of varieties that does not seem to have been studied yet, but for
which the existing technology is likely to prove useful.

Let us consider the type \texttt{xv}  surface 
$$
S_1=\{[x_0,\ldots,x_4]\in \bfP^4: x_0x_1+x_2^2=x_0x_4+x_1x_2+x_3^2 = 0\},
$$
in more detail.  Now we have already seen that \eqref{c2} holds for $S_1$.
Nonetheless, la Bret\`eche and Browning \cite{math.NT/0412086} have made an
exhaustive study of $S_1$, partly in an attempt to 
lay down a template for the treatment of other surfaces in the
table.  In doing so several new features have been  revealed.  
For $s\in \C$ such that $\Re e(s)>1$, let 
\begin{equation}\lab{height-zeta}
Z_{U,H}(s):=\sum_{x \in U(\Q)}H(x)^{-s}
\end{equation}
denote the corresponding height zeta function, 
where $U=U_1$ denotes the open subset formed by deleting the unique line
$x_0=x_2=x_3=0$ from $S_1$. 
The analytic properties of $Z_{U_1,H}(s)$ are intimately related to
the asymptotic behaviour of the counting function $N_{U_1,H}(B)$, and
it is relatively straightforward to translate between them.  For
$\sigma \in \R$, let $\HH_\sigma$ denote the half-plane $\{s \in \C:
\Re e(s)>\sigma\}$. Then with this notation in mind we have the following result \cite[Theorem
1]{math.NT/0412086}. 

\begin{thm}\lab{4-1}
There exists a constant $\alpha \in \R$, a function $F_1(s)$ that is
meromorphic on $\HH_{9/10}$ with a pole of order $6$ at $s=1$, and a function
$F_2(s)$ that is holomorphic on $\HH_{5/6}$, such that 
$$
Z_{U_1,H}(s) =F_1(s) +\alpha (s-1)^{-1}+F_2(s),
$$
for $s \in \HH_1$.  In particular $Z_{U_1,H}(s)$ has an analytic continuation to $\HH_{9/10}.$ 
\end{thm}

It should be highlighted  that there exist remarkably precise descriptions of
the terms $F_1,F_2,\al$ that appear in the statement of the theorem.
An application of Perron's formula enables one to deduce a
corresponding asymptotic formula for $N_{U_1,H}(B)$ that verifies
\eqref{c2}, with $\rho_{S_1}=6$.  Actually one is led to the much
stronger statement that there exists a polynomial $f$ of degree $5$ such that 
for any $\delta \in (0,1/12)$ we have
\begin{equation}\lab{c3}
N_{U,H}(B)=B f(\log B) +O(B^{1-\delta}),
\end{equation}
with $U=U_1$, in which the leading coefficient of $f$ agrees
with Peyre's prediction.

No explicit use is made of the fact that $S_1$ is an equivariant compactification of
$\mathbb{G}_a^2$ in the proof of Theorem \ref{4-1}, and this renders
the method applicable to other surfaces in the list that are not of
this type.  For example, 
in further work la Bret\`eche and Browning \cite{math.NT/0502510}
have also established the Manin conjecture for the $\Dfour$ surface 
$$
S_2=\{[x_0,\ldots,x_4]\in \bfP^4: x_0x_1+x_2^2=x_0^2+x_1x_4+x_3^2 = 0\},
$$
which represents the type \texttt{xiii} surface in the table. 
This surface is not split, since it contains 
the pair of lines $x_1=x_2=x_0\pm  i x_3=0$, and it turns out that 
$\Pic(\tS_2)$ has rank $4$.  In fact $\tS_2$ has singularity 
type $\mathbf{C}_3$ over $\Q$, in the sense of Lipman \cite[\S
24]{lipman}, which becomes a $\Dfour$ singularity over $\overline{\Q}$.
Building on the techniques developed in the proof of Theorem
\ref{4-1}, a result of the same quality is obtained for the
corresponding height zeta function $Z_{U_2,H}(s)$, and this leads to
an estimate of the shape \eqref{c3} for $\delta \in (0,3/32)$, with $U=U_2$
and $\deg f=3$.

One of the aims of this survey is to give an overview of the
various ideas and techniques that have been used to study the surfaces
$S_1,S_2$ above.  We shall illustrate the basic method by
giving a simplified analysis of a new example from the table.  Let us consider
the $3\textbf{A}_1$ surface 
\begin{equation}\lab{S-new}
S_3=\{[x_0,\ldots,x_4]\in \bfP^4: x_0x_1+x_2^2=x_1x_2+x_2^2 +x_3x_4= 0\},
\end{equation}
which represents the type \texttt{v} surface in the table, and is
neither toric, nor
an equivariant compactification of $\mathbb{G}_a^2$.  The surface
has singularities at the points $[1,0,0,0,0], [0,0,0,1,0]$ and $
[0,0,0,0,1]$, and contains precisely $6$ lines
$$
x_i=x_2=x_j=0, \quad x_0+x_2=x_1+x_2=x_j=0,
$$
where $i \in \{0,1\}$ and $j \in \{3,4\}$. Since $S_3$ is split, one finds that 
the expected exponent of $\log B$ in \eqref{c2} is $\rho_{S_3}-1=5$. 
We shall establish the following result.

\begin{thm}\lab{4-2} We have
$
N_{U_3,H}(B)=O\big(B (\log B)^5\big).
$
\end{thm}

As pointed out to the author by R\'egis de la Bret\`eche, it is
possible to establish a corresponding lower bound 
$N_{U_3,H}(B)\gg B (\log B)^5$, using little more than the most basic
estimates for integers restricted to lie in fixed congruence classes.
In fact, with
more work, it ought even to be possible to obtain an asymptotic formula for 
$N_{U_3,H}(B)$.
In the interests of brevity, however, we have chosen to pursue neither
of these problems here.

\subsection{Del Pezzo surfaces of degree $3$}\lab{1.3}

The del Pezzo surfaces $S\subset \bfP^3$ of degree $3$ are readily recognised as
the geometrically integral cubic surfaces in $\bfP^3$, that are not cones.  Given such a
surface $S$ defined over $\Q$, we may always find an absolutely
irreducible cubic form $C(\x)\in \Z[x_0,x_1,x_2,x_3]$ such that $S=\Proj(\Q[\x]/(C))$.
Let us begin by considering the situation for non-singular cubic
surfaces. 
In this setting $U\subset S$ is taken to be the open subset formed by
deleting the famous $27$ lines from $S$.  Although
Peyre and Tschinkel \cite{p-t1, p-t2} have provided 
ample numerical evidence for the validity of the Manin conjecture for
non-singular cubic surfaces,  we are unfortunately still rather far away from
proving it.  The best upper bound available is $\NU(B)=O_{\varepsilon,S}(B^{4/3+\varepsilon})$,
due to Heath-Brown \cite{MR98h:11083}.  This applies when the surface
$S$ contains $3$ coplanar lines defined over $\Q$, and in particular
to the {\em Fermat cubic surface} $x_0^3+x_1^3=x_2^3+x_3^3$.
The problem of proving lower bounds is somewhat easier.  Under the
assumption that $S$ contains a pair of skew lines defined over $\Q$,
Slater and Swinnerton-Dyer \cite{s-swd} have shown that $\NU(B)\gg_S
B(\log B)^{\rho_S-1}$, as predicted by the Manin conjecture.
This does not apply to the Fermat cubic surface, however,
since the only skew lines contained in this surface are defined over
$\Q(\sqrt{-3})$.  It would be interesting to extend the work of Slater
and Swinnerton-Dyer to cover such cases.

Much as in the previous section, it turns out that far better
estimates are available for singular cubic surfaces.  The
classification of such surfaces is a well-established subject, and essentially
goes back to the work of Cayley \cite{cayley} and Schl\"afli
\cite{cayley'} over a century ago.   
A contemporary classification of singular cubic surfaces, using the
terminology of modern classification theory, has since been given
by Bruce and Wall \cite{MR80f:14021}.
As in the previous section, the Manin conjecture is already
known to hold for several of these surfaces by virtue of the fact that
they are equivariant compactifications of $\mathbb{G}_a^2$, or toric,
such as the $3\textbf{A}_2$ surface 
\begin{equation}\lab{3a2}
S_4=\{[x_0,x_1,x_2,x_3]\in \bfP^3: x_0^3=x_1x_2x_3\},
\end{equation}
for example. In fact a number of authors have studied this
particular surface, including la Bret\`{e}che \cite{MR2000b:11074}, 
Fouvry  \cite{MR2000b:11075}, and Heath-Brown and Moroz \cite{MR2000f:11080}.
Of the asymptotic formulae obtained, the most impressive is the first, which consists of an
estimate like \eqref{c3} for any $\delta \in (0,1/8)$, with
$U=U_4\subset S_4$ and a suitable polynomial $f$ of degree $6$.  The next surface to
receive serious attention was the so-called {\em Cayley cubic surface}
\begin{equation}\lab{cayley-eq}
S_5=\{[x_0,x_1,x_2,x_3]\in \bfP^3: x_0x_1x_2+x_0x_1x_3+x_0x_2x_3+x_1x_2x_3=0\},
\end{equation}
of singularity type $4\mathbf{A}_1$.  This contains $9$ lines, all of
which are defined over $\Q$, and Heath-Brown \cite{0981.32025} has shown that
there exist absolute constants $A_1,A_2>0$ such that 
$$
A_1 B (\log B)^6 \leq N_{U_5,H}(B) \leq A_2 B (\log B)^6,
$$
where $U_5\subset S_5$ is the usual open subset.  An estimate of
precisely the same form has also been obtained by Browning
\cite{math.NT/0404245} for the $\Dfour$ surface
$$
S_6=\{[x_0,x_1,x_2,x_3]\in \bfP^3:  x_1x_2x_3=x_0(x_1+x_2+x_3)^2\}.
$$
In both cases the corresponding Picard group has rank $7$, so that the
exponents of $B$ and $\log B$ agree with Manin's
prediction.

The final surface to have been studied extensively is the $\Esix$ cubic surface
\begin{equation}\lab{e6}
S_7=\{[x_0,x_1,x_2,x_3] \in \bfP^3: x_1x_2^2+x_2x_0^2+x_3^3 = 0\},
\end{equation}
which contains a unique line $x_2=x_3 = 0$.  
Let $U_7\subset S_7$ denote the open subset formed by deleting
the line from $S_7$, and recall the notation \eqref{height-zeta} 
for the height zeta function $Z_{U_7,H}(s)$ and that of the 
half-plane $\HH_\sigma$ introduced before Theorem~\ref{4-1}. Then recent work of la Bret\`eche, Browning  and
Derenthal \cite{e6} has succeeded in establishing the following result.

\begin{thm}\lab{3-1}
There exists a constant $\alpha \in \R$, a function $F_1(s)$ that is
meromorphic on $\HH_{9/10}$ with a pole of order $7$ at $s=1$, and a function
$F_2(s)$ that is holomorphic on $\HH_{43/48}$, such that 
$$
Z_{U_7,H}(s) =F_1(s) +\alpha (s-1)^{-1}+F_2(s),
$$
for $s \in \HH_1$.  In particular $Z_{U_7,H}(s)$ has an analytic continuation to $\HH_{9/10}.$ 
\end{thm}

As in Theorem \ref{4-1}, the terms $F_1,F_2,\al$ have a very explicit description.
An application of Perron's formula now yields an
asymptotic formula of the shape \eqref{c3} for $\delta \in (0,1/11)$,
with $U=U_7$  and a suitable polynomial $f$ of degree $6$.
This too is in complete agreement with the Manin conjecture.
It should be remarked that, in work to appear, Michael Joyce has 
independently established the Manin conjecture 
for $S_7$ in his doctoral thesis at Brown University, albeit only with
a weaker error term of  $O(B (\log B)^5)$.

\section{Refinements of the Manin conjecture}

The purpose of this section is to consider in what way one might hope
to refine the conjecture of Manin.  We have already seen a number of
examples in which asymptotic formulae of the shape \eqref{c3} hold, and
it is very natural to suppose that this is the case for any (possibly
singular) del Pezzo surface $S \subset \bfP^d$ of degree $d$, where
as usual $U\subseteq S$ denotes the open subset formed by deleting any
exceptional divisors from $S$, and $\rho_S$ denotes the rank of the
Picard group of $S$ (possibly of $\tS$). Let us record this formally here.

\begin{con2}
Let $S,U, \rho_S$ be as above.  Then there exists $\delta>0$, and a
polynomial $f\in \R[x]$ of degree $\rho_S-1$,  such that \eqref{c3} holds.
\end{con2}

The leading coefficient of $f$ should of course agree with the
prediction of Peyre et al. It would be interesting to gain a
conjectural understanding of the lower order coefficients of $f$,
possibly in terms of the geometry of $S$.  
At this stage it seems worth drawing attention to the surprising nature
of the constants $\alpha$ that appear in Theorems \ref{4-1} and
\ref{3-1}, not least because they contribute to the constant coefficient of 
$f$.  In both cases we have $\alpha=\frac{12}{\pi^2}+\beta$, where the first term corresponds to
an isolated conic in the surface, and the second is
purely arithmetic in nature and takes a very complicated shape (see
\cite[Eq. (5.25)]{math.NT/0412086} and \cite[Eq. (8.49)]{e6}). 
It arises through the error in  
approximating certain arithmetic quantities by real-valued
continuous functions, and involves the application of results about
the equidistribution of squares in fixed residue classes. 

One might also ask what one expects to be the true order of magnitude of
the error term in \eqref{c3}.  This a question that Swinnerton-Dyer
has recently addressed \cite[Conjecture $2$]{swd}, inspired by
comparisons with the explicit formulae from prime number theory.

\begin{con3}
Let $S,U, \rho_S$ be as above.  Then there exist positive contstants 
$\theta_1,\theta_2,\theta_3<1$ with $\theta_1<\min \{\theta_2,\theta_3\}$, a
polynomial $f\in \R[x]$ of degree $\rho_S-1$,  a constant $\gamma \in
\R$, and a sequence of $\gamma_n\in \C$, such for any $\ve>0$ we have
$$
N_{U,H}(B)=B f(\log B) +\gamma B^{\theta_3}+\Re e \sum \gamma_n B^{\theta_2+it_n}+O_\ve(B^{\theta_1+\ve}).
$$
Here $\frac{1}{2}+it_n$ runs through a set of non-trivial zeros of the Riemann
zeta function, with the $t_n$ positive and monotonic increasing,
such that $\sum |\gamma_n|^2$ and $\sum t_n^{-2}$ are convergent.
\end{con3}

In fact Swinnerton-Dyer formulates the conjecture for non-singular
cubic surfaces, with $\theta_1<\frac{1}{2}=\theta_2$ and $\gamma=0$. 
There is no reason, however, to expect that it doesn't hold more
generally, and one might even suppose that the constants $\theta_2,\theta_3$
somehow relate to the nature of the surface singularities.  In work
currently under preparation, la Bret\`eche and Swinnerton-Dyer have
provided significant evidence for this finer conjecture for the singular
cubic surface \eqref{3a2}. Under the assumption of the 
Riemann hypothesis it is shown that the conjecture holds for $S_4$,
with $(\theta_1,\theta_2,\theta_3)=(\frac{4}{5},\frac{13}{16},
\frac{9}{11})$ and $\gamma \neq 0$.

\section{Available tools}

There are a variety of tools that can be brought to bear upon the problem of estimating
the counting function \eqref{count} for appropriate subsets $U$ of
projective algebraic varieties.  Most of these are rooted in analytic number theory.
When the dimension of the variety is
large compared to its degree, the Hardy--Littlewood circle method can
often be applied successfully (see Davenport \cite{dav}, for
example).  When the variety has a suitable ``cellular'' structure,
techniques involving harmonic analysis on adelic groups can be
employed (see Tschinkel \cite{t-2002}, for example).
We shall say nothing more about these methods here, save to observe
that outside of the surfaces covered by the collective work of
Batyrev, Chambert-Loir and Tschinkel \cite{MR1620682, ct}, 
they do not seem capable of establishing the Manin conjecture for all
del Pezzo surfaces.   

In fact we still have no clear vision of which methods are most
appropriate, and it is conceivable that the methods needed to handle
the singular del Pezzo surfaces of low degree are quite different from 
those needed to handle the non-singular surfaces. Given our inability
to prove the Manin conjecture for a single non-singular del Pezzo
surface of degree $3$ or $4$, we shall say no more about them here,
save to observe that the sharpest results we have are for examples
containing conic bundle structures over the ground field.  
Instead we shall concentrate on the situation for singular del Pezzo surfaces of degree $3$ or $4$.
Disappointing as it may seem, it is hard to imagine that we will see
how to prove Manin's conjecture for all del Pezzo
surfaces without first attempting to do so for a number of very
concrete representative examples.  
As a cursory analysis of the proofs of 
Theorems \ref{4-1}--\ref{3-1} shows, the techniques that have been successfully applied so far
are decidedly ad-hoc.  Nonetheless there are a few salient features that are worthy of
amplification, and this will be the focus of the two subsequent sections.

\subsection{The universal torsor}

Universal torsors were originally introduced 
by Colliot-Th\'el\`ene and Sansuc \cite{ct1,ct2} 
to aid in the study of the Hasse principle and weak approximation for
rational varieties.  Since their inception it is now well-recognised that
they also have a central r\^ole to play in proofs of the Manin conjecture for Fano varieties.
Let $S \subset \bfP^d$ be a del Pezzo
surface of degree $d\in \{3,4,5\}$, and let $\tS$ denote the minimal desingularisation
of $S$ if it is singular, and $\tS=S$ otherwise.  Let $E_1,\ldots,E_{10-d} \in
\rom{Div} (\tS)$ be generators for the geometric Picard group of $\tS$, and
let $E_i^\times=E_i\setminus\{\mbox{zero section}\}$.  
Working over $\overline{\Q}$,  ``the'' {\em universal torsor} of $\tS$ is given by the
action of $\mathbb{G}_m^{10-d}$ on the map
$$
\pi: E_1^\times \times_\tS \cdots \times_\tS E_{10-d}^\times
\rightarrow \tS.
$$
In practice this action can be made completely explicit, thereby giving equations for the
universal torsor.  A proper discussion of universal torsors would
take us too far afield at present, and the reader should
consult the survey of Peyre \cite{MR2029862} for further details, or indeed the
construction of Hassett and Tschinkel \cite{MR2029868}. The latter outlines
an alternative approach to universal torsors via the Cox ring.
The guiding principle behind the use of universal torsors is simply that they ought to be
arithmetically simpler than the original variety.
The universal torsors that we shall encounter all have embeddings as
open subsets of affine varieties of higher dimension,  and the
general theory ensures that there is a bijection between $U(\Q)$
--- where $U \subset S$ is the usual open subset formed by deleting the lines
from $S$ --- and a suitable set of integral points on the
corresponding universal torsor. 
We shall see shortly how one may often use arguments from elementary
number theory to explicitly derive these bijections.

Let us begin by giving a few examples.  In the proof of Theorem
\ref{5-1} a passage to the universal torsor is a crucial first step,
and was originally carried out by Salberger in his unpublished proof of
the bound $N_{U_0,H}(B)=O(B(\log B)^4)$, announced in the 
{\em Borel seminar} at Bern in 1993.
Recall the Pl\"ucker embedding 
$$
z_{i,j}z_{k,\ell}-z_{i,k}z_{j,\ell}+ z_{i,\ell}z_{j,k}=0, 
$$
where $\{i,j,k,\ell\}$ runs through cyclic permutations of $\{1,2,3,4,5\}$, 
of the Grassmannian $Gr(2,5)\subset \bfP^9$ of $2$-dimensional linear subspaces of $\Q^5.$
Then the universal torsor $\pi: \mcal{T}_0\rightarrow S_0$ above $S_0$ is a certain open subset of the
affine cone over $Gr(2,5)$. To count points of bounded height in
$U_0(\Q)$ it is then enough to count integral points $(z_{i,j})_{1\leq i<j \leq 5} \in {\Z_*}^9$ 
on this cone, where $\Z_*:=\Z \setminus \{0\}$, subject to a number
of side conditions.  A thorough account of
this particular example, and how it extends to arbitrary del Pezzo
surfaces of degree $5$ can be found in the work of Skorobogatov \cite{skor}.
A second example is calculated by Hassett and Tschinkel
\cite{MR2029868} for the $\textbf{E}_6$ cubic surface \eqref{e6}.
There it is shown that the universal torsor above $\tS_7$ has the equation
\begin{equation}\lab{ut-e6}
\tau_\ell\xi_\ell^3\xi_4^2\xi_5+\tau_2^2\xi_2+\tau_1^3\xi_1^2\xi_3 = 0,
\end{equation}
for variables
$\tau_1,\tau_2,\tau_\ell,\xi_1,\xi_2,\xi_3,\xi_\ell,\xi_4,\xi_5,\xi_6$.
One of the variables does not explicitly appear in \eqref{ut-e6}, and
the torsor should be thought of as being embedded in $\A^{10}$.
The universal torsors that turn up in the proofs of Theorems
\ref{4-1} and \ref{4-2} can also be embedded in affine space via a single
equation.

We proceed to carry out explicitly the passage to the universal torsor
for the $3\textbf{A}_1$ surface  \eqref{S-new}.
We shall use $\N$ to denote the set of positive integers, and for any 
$n \geq 2$ we let $Z^{n}$ denote the set of {\em primitive} vectors in $\Z^{n}$,
by which we mean that the greatest common divisor of the components should be $1$.
We may clearly assume that $S_3$ is defined by the forms
$Q_1(\x)=x_0x_1-x_2^2$ and $Q_2(\x)=x_2^2-x_1x_2+x_3x_4.$ 
Now if $x\in U_3(\Q)$ is represented by the vector
$\x \in Z^5$, then $x_0\cdots x_4\neq 0$ and
$H(x)=\max\{|x_0|,|x_1|,|x_3|,|x_4|\}.$
Moreover, $x_0$ and $x_1$ must share the same sign. On
taking $x_0,x_1$ to both be positive, and noting that 
$\x$ and $-\x$ represent the same point in $\bfP^4$, we deduce that
$$
N_{U_3,H}(B)= \#\big\{\x\in Z^5:  0<x_0,x_1,|x_3|,|x_4| \leq B, ~ 
Q_1(\x)=Q_2(\x)=0\big\}.
$$
Let us begin by considering solutions $\x \in Z^5$ to the equation 
$Q_1(\x)=0.$ There is a bijection
between the set of integers $x_0,x_1,x_2$ such that
$x_0,x_1>0$ and $x_0x_1=x_2^2$, and the set of $x_0,x_1,x_2$ such that
$x_0=z_0^2z_2, x_1=z_1^2z_2$ and $x_2=z_0z_1z_2,$
for non-zero integers $z_0, z_1,z_2$ such that $z_1,z_2>0$ and
$\hcf(z_0,z_1)=1.$
We now substitute these values into the equation
$Q_2(\x)=0,$ 
in order to obtain 
\begin{equation}\lab{golf-1}
z_0^2z_1^2z_2^2-z_0z_1^3z_2^2+x_3x_4=0.
\end{equation}
It follows from the coprimality relation $\hcf(x_0,\ldots,x_4)=1$ that we
also have $\hcf(z_2,x_3,x_4)=1$.
Now we may conclude from \eqref{golf-1} that $z_0z_1^2z_2^2$ divides $x_3x_4$.
Let us write $y_1=\hcf(z_1,x_3,x_4)$ and 
$z_1=y_1y_1', x_3=y_1y_3', x_4=y_1y_4',$
with  $y_1,y_1',y_3',y_4'$ non-zero integers such that $y_1,y_1'>0$ and 
$\hcf(y_1',y_3',y_4')=1$.
Then $z_0y_1'^2z_2^2$ divides $y_3'y_4'$.  We now write
$z_0=y_{03}y_{04}, y_3'=y_{03}y_3$ and $y_4'=y_{04}y_4,$
for non-zero integers $y_{03},y_{04},y_3,y_4$. We therefore conclude that $y_1'^2z_2^2$ divides
$y_3y_4$, whence 
there exist positive integers $y_{13},y_{14},y_{23},y_{24}$ and non-zero integers
$y_{33},y_{34}$ such that
$y_1'=y_{13}y_{14}, z_2=y_{23}y_{24}, 
y_3=y_{13}^2y_{23}^2y_{33}$ and $y_4=y_{14}^2y_{24}^2y_{34}.$
Substituting these into \eqref{golf-1} yields the equation
\begin{equation}\lab{ut}
y_{03}y_{04}-y_1y_{13}y_{14}+y_{33}y_{34}=0.
\end{equation}
This equation gives an affine embedding of the universal torsor over
$\tS_3$, though we shall not prove it here.
Furthermore, we may combine all of the various coprimality relations
above to deduce that
\begin{equation}\lab{polo-3}
\hcf(y_{13}y_{14}y_{23}y_{24},y_{13}y_{23}y_{33},y_{14}y_{24}y_{34})=1,
\end{equation}
and
\begin{equation}\lab{polo-4}
\hcf(y_{03}y_{04},y_{13}y_{14})=\hcf(y_{1},y_{03}y_{04}y_{23}y_{24})= 1.
\end{equation}

At this point we may summarize our argument as follows.  
Let $\mcal{T}$ denote the set of non-zero integer vectors
$\ma{y}=(y_1,y_{03},y_{04},y_{13},y_{14},y_{23},y_{24},y_{33},y_{34})$
such that \eqref{ut}--\eqref{polo-4}
all hold, with $y_1, y_{13},y_{14},y_{23},y_{24}>0$.
Then for any 
$\x \in Z^5$ such that $Q_1(\x)=Q_2(\x)=0$ and $x_0,x_1,|x_3|,|x_4|>0$,
we have shown that there exists
$\ma{y}\in \mcal{T}$ such that 
$$
\begin{array}{l}
x_0 = y_{03}^2y_{04}^2y_{23}y_{24},\\
x_1 = y_1^2y_{13}^2y_{14}^2y_{23}y_{24},\\ 
x_2 = y_1y_{03}y_{04}y_{13}y_{14}y_{23}y_{24}, \\
x_3 = y_1y_{03}y_{13}^2y_{23}^2y_{33}, \\
x_4 = y_1y_{04}y_{14}^2y_{24}^2y_{34}.
\end{array}
$$
Conversely, it is not hard to check that given any
$\ma{y} \in \mcal{T}$ the point $\x$ given above
will be a solution of the equations $Q_1(\x)=Q_2(\x)=0$, with
$\x \in Z^5$ and $x_0,x_1,|x_3|,|x_4|>0$.  
Let us define the function $\Psi: \R^9 \rightarrow \R_{\geq 0}$, given by
$$
\Psi(\ma{y})=
\max\left\{
\begin{array}{l}
|y_{03}^2y_{04}^2y_{23}y_{24}|, \quad 
|y_1^2y_{13}^2y_{14}^2y_{23}y_{24}|,\\
|y_1y_{03}y_{13}^2y_{23}^2y_{33}|, \quad 
|y_1y_{04}y_{14}^2y_{24}^2y_{34}|
\end{array}
\right\}.
$$
Then we have established the following result.

\begin{lem}\lab{base}
We have
$
N_{U_3,H}(B)=\#\big\{\ma{y} \in \mcal{T}:
     \Psi(\ma{y}) \leq B\big\}.
$
\end{lem}

In this section we have given several examples of universal torsors,
and we have ended by demonstrating how elementary number theory can
sometimes be used to calculate them with very little trouble.  In fact 
the general machinery of Colliot-Th\'el\`ene--Sansuc 
\cite{ct1, ct2}, or that of Hassett--Tschinkel \cite{MR2029868},
essentially provides an algorithm for calculating the universal torsor
over any singular del Pezzo surface of degree $3$ or $4$.
It should be stressed, however, that if this 
constitutes being given the keys to the city, it does not tell us
where in the city the proof is hidden.

\subsection{The next step}\lab{3.2}

The purpose of this section is to overview some of the techniques that
have been developed for counting integral points on the
parametrization that arises out of the passage to the universal
torsor, as discussed above.  In the proofs of Theorems
\ref{5-1}--\ref{3-1} the torsor equations all take the shape
$$
A_j+B_j+C_j=0, \quad (1\leq j \leq J),
$$
for monomials $A_j,B_j,C_j$ of various degrees in the
appropriate variables.   By fixing some of the variables at
the outset, one is 
then left with the problem of counting integer solutions to a system of
Diophantine equations, subject to certain constraints.  If one is
sufficiently clever about which variables to fix first, then one can
sometimes be left with a quantity 
that we know how to estimate ---  and crucially --- for which we can control
the overall contribution from the error term when it is summed over
the remaining variables. 

Let us sketch this phenomenon briefly with the torsor equation
\eqref{ut-e6} that is used in the proof of Theorem \ref{3-1}.
It turns out that the way to proceed here is to fix all of the
variables apart from $\tau_1,\tau_2,\tau_\ell$.
One may then  view the  equation as a congruence
$$
\tau_2^2\xi_2\equiv -\tau_1^3\xi_1^2\xi_3 ~~\mod{\xi_\ell^3\xi_4^2\xi_5},
$$
in order to take care of the summation over $\tau_\ell$.
This allows us to employ very standard facts about the number of integer solutions to
polynomial congruences that are restricted to lie in certain regions,
and this procedure yields a main term and an error term which the
remaining variables need to be summed over.  However, while the treatment of
the main term is relatively routine, the treatment of the error term presents a much
more serious obstacle.  Although we do not have space to discuss it in
any detail, it is here that the unexpected constant $\alpha$ arises in
Theorem \ref{3-1} (and, indeed, in Theorem \ref{4-1}).

The sort of approach discussed above, and more generally the
application of lattice methods to count solutions to ternary
equations,  is a very useful one.  It plays a crucial role in the proof of
the following result due to Heath-Brown \cite[Lemma 3]{cayley},  which 
forms the next ingredient in our proof of Theorem \ref{4-2}.

\begin{lem}\lab{cayley-teq}
Let $K_1,\ldots,K_7\geq 1$ be given, and let $\mcal{M}$ denote the
number of non-zero solutions $m_1,\ldots,m_7\in\Z$ to the equation
\[
m_1m_2-m_3m_4m_5+ m_6m_7=0,
\]
subject to the conditions
$K_k<|m_k|\leq 2K_k$ for $1\leq k\leq 7$, and 
\begin{equation}\label{4.1}
\hcf(m_1m_2, m_3m_4m_5)=1.
\end{equation}
Then we have $\mcal{M}\ll K_1K_2K_3K_4K_5$.
\end{lem}

For comparison, we note that it is a trivial matter to establish the 
bound $\mcal{M}\ll_\ve (K_1K_2K_3K_4K_5)^{1+\ve}$,
using standard estimates for the divisor function.

\subsection{Completion of the proof of Theorem \ref{4-2}}

We are now ready to complete the proof of Theorem \ref{4-2}.
We shall begin by estimating the contribution to $N_{U_3,H}(B)$ from the values of $\ma{y}$
appearing in Lemma~\ref{base} that are constrained to lie in a certain
region.  Let $Y_1,Y_{i3},Y_{i4}\geq 1$, where throughout this section
$i$ denotes a generic index from the set $\{0,1,2,3\}$.  Then we write
$\mcal{N}=\mcal{N}(Y_1,Y_{03},Y_{04},Y_{13},Y_{14},Y_{23},Y_{24},Y_{33},Y_{34})$
for the total contribution to $N_{U_3,H}(B)$ from $\ma{y}$ satisfying 
\begin{equation}\lab{range1}
Y_1\leq y_1 <2Y_1, \quad 
Y_{i3}\leq |y_{i3}|<2Y_{i3}, \quad 
Y_{i4}\leq |y_{i4}|<2Y_{i4}.
\end{equation}
Clearly it follows from the inequality $\Psi(\ma{y}) \leq B$ that 
$\mcal{N}=0$ unless
\begin{equation}\lab{r1}
Y_{03}^2Y_{04}^2Y_{23}Y_{24} \ll B, \quad 
Y_1^2Y_{13}^2Y_{14}^2Y_{23}Y_{24} \ll B,
\end{equation}
and 
\begin{equation}\lab{r2}
Y_1Y_{03}Y_{13}^2Y_{23}^2Y_{33}\ll B, \quad 
Y_1Y_{04}Y_{14}^2Y_{24}^2Y_{34}\ll B.
\end{equation}
In our estimation of $N_{U_3,H}(B)$, we may clearly assume without
loss of generality that 
\begin{equation}\lab{mendip}
Y_{03}Y_{13}^2Y_{23}^2Y_{33}\leq Y_{04}Y_{14}^2Y_{24}^2Y_{34}.
\end{equation}

We proceed to show how the equation \eqref{ut} forces certain
constraints upon the choice of dyadic ranges in
\eqref{range1}.  There are three basic cases that can occur.  Suppose first that 
\begin{equation}\lab{case-1a}
c_2Y_{03}Y_{04} \leq Y_1Y_{13}Y_{14},
\end{equation} 
for an absolute constant $c_2>0$.
Then it follows from \eqref{ut} that 
\begin{equation}\lab{case1}
Y_{33}Y_{34} \ll Y_{1}Y_{13}Y_{14} \ll Y_{33}Y_{34},
\end{equation}
provided that $c_2$ is chosen to be sufficiently large.
Next, we suppose that
\begin{equation}\lab{case-2a}
c_1Y_{03}Y_{04} \geq Y_1Y_{13}Y_{14},
\end{equation} 
for an absolute constant $c_1>0$. Then we may deduce from \eqref{ut} that 
\begin{equation}\lab{case2}
Y_{33}Y_{34} \ll Y_{03}Y_{04} \ll Y_{33}Y_{34},
\end{equation}
provided that $c_1$ is chosen to be sufficiently small.
Let us henceforth assume that the values of $c_1,c_2$ are fixed in
such a way that \eqref{case1} holds, if \eqref{case-1a} holds, and
\eqref{case2} holds, if \eqref{case-2a} holds.
Finally we are left with the possibility that 
\begin{equation}\lab{case3}
c_1Y_{03}Y_{04} \leq Y_1Y_{13}Y_{14} \leq c_2Y_{03}Y_{04}.
\end{equation}
We shall need to treat the cases \eqref{case-1a}, \eqref{case-2a} and
\eqref{case3} separately.  

We shall take $\ma{m}_{j,k}=(y_1,y_{13},y_{14},y_{j3},y_{j4},y_{k3},y_{k4})$
in our application of Lemma \ref{cayley-teq}, for $(j,k)=(0,3)$ and $(3,0)$.
In particular the
coprimality relation \eqref{4.1} follows directly from \eqref{polo-3} and
\eqref{polo-4}, and we may conclude that 
\begin{equation}\lab{N1}
\mcal{N} \ll Y_1Y_{13}Y_{14}Y_{23}Y_{24} \min\{Y_{03}Y_{04},Y_{33}Y_{34}\},
\end{equation}
on summing over all of the available $y_{23},y_{24}$.
It remains to sum this contribution over the various dyadic intervals
$Y_1,Y_{i3},Y_{i4}$.  Suppose for the moment that we are interested in summing over all
possible dyadic intervals $X \leq |x|<2X$, for which $|x| \leq
\mcal{X}$.  Then there are plainly $O(\log \mcal{X})$ possible choices
for $X$.  In addition to this basic estimate, we shall make frequent
use of the estimate $\sum_{X} X^{\delta} \ll_\delta \mcal{X}^{\delta},$
for any $\delta>0$.

We begin by assuming that \eqref{case-1a} holds, so that \eqref{case1} also holds. 
Then we may combine \eqref{mendip} with \eqref{case1} in order to deduce that
$$
Y_{13} \ll \min\left\{  \frac{Y_{04}^{1/2}Y_{14}Y_{24}Y_{34}^{1/2}}{Y_{03}^{1/2}Y_{23}Y_{33}^{1/2}}, 
\frac{Y_{33}Y_{34}}{Y_1Y_{14}}\right\} \ll
\frac{Y_{04}^{1/4}Y_{24}^{1/2}Y_{33}^{1/4}Y_{34}^{3/4}}{Y_1^{1/2}Y_{03}^{1/4}Y_{23}^{1/2}}.
$$
We may now apply \eqref{N1} to obtain 
\begin{align*}
\sum_{\colt{Y_1,Y_{i3},Y_{i4}}{\mbox{\scriptsize{\eqref{case-1a} holds}}}} \mcal{N} 
&\ll 
\sum_{\colt{Y_1,Y_{i3},Y_{i4}}{\mbox{\scriptsize{\eqref{case-1a} holds}}}}
Y_1Y_{03}Y_{04}Y_{13}Y_{14}Y_{23}Y_{24}\\
&\ll \sum_{\colt{Y_{03},Y_{04},Y_{33},Y_{34}}{Y_1,Y_{14},Y_{23},Y_{24}}} 
Y_1^{1/2}Y_{03}^{3/4}Y_{04}^{5/4}Y_{14}Y_{23}^{1/2}Y_{24}^{3/2}Y_{33}^{1/4}Y_{34}^{3/4}.
\end{align*}
But now \eqref{r2} implies that 
$Y_{14}\ll B^{1/2}/(Y_1^{1/2}Y_{04}^{1/2}Y_{24}Y_{34}^{1/2})$, 
and \eqref{case1} and \eqref{case-1a} together imply that $Y_{03}\ll
Y_{33}Y_{34}/Y_{04}$.  We therefore deduce that
\begin{align*}
\sum_{\colt{Y_1,Y_{i3},Y_{i4}}{\mbox{\scriptsize{\eqref{case-1a} holds}}}} \mcal{N} 
&\ll B^{1/2}\sum_{\colt{Y_{03},Y_{04},Y_{33}}{Y_1,Y_{23},Y_{24},Y_{34}}} 
Y_{03}^{3/4}Y_{04}^{3/4}Y_{23}^{1/2}Y_{24}^{1/2}Y_{33}^{1/4}Y_{34}^{1/4}\\
&\ll B^{1/2}\sum_{\colt{Y_1,Y_{04},Y_{33}}{Y_{23},Y_{24},Y_{34}}} 
Y_{23}^{1/2}Y_{24}^{1/2}Y_{33}Y_{34}.
\end{align*}
Finally it follows from \eqref{r1} and \eqref{case1} that $Y_{33}\ll
B^{1/2}/(Y_{23}^{1/2}Y_{24}^{1/2}Y_{34})$, whence
\begin{align*}
\sum_{\colt{Y_1,Y_{i3},Y_{i4}}{\mbox{\scriptsize{\eqref{case-1a} holds}}}} \mcal{N} 
\ll B \sum_{Y_{04},Y_{13},Y_{14},Y_{23},Y_{34}} 1 \ll B (\log B)^5,
\end{align*}
which is satisfactory for the theorem.  

Next we suppose that \eqref{case-2a} holds,  so that \eqref{case2} also holds. 
In this case it follows from \eqref{mendip}, together with the inequality $Y_1Y_{13}Y_{14}\ll Y_{03}Y_{04}$,
that
$$
Y_{13} \ll \min\left\{  \frac{Y_{04}^{1/2}Y_{14}Y_{24}Y_{34}^{1/2}}{Y_{03}^{1/2}Y_{23}Y_{33}^{1/2}}, 
\frac{Y_{03}Y_{04}}{Y_1Y_{14}}\right\} \ll
\frac{Y_{03}^{1/4}Y_{04}^{3/4}Y_{24}^{1/2}Y_{34}^{1/4}}{Y_1^{1/2}Y_{23}^{1/2}Y_{33}^{1/4}}.
$$
On combining this with the inequality
$Y_{14}\ll B^{1/2}/(Y_1^{1/2}Y_{04}^{1/2}Y_{24}Y_{34}^{1/2})$, that
follows from \eqref{r2}, we may therefore deduce from \eqref{N1} that
\begin{align*}
\sum_{\colt{Y_1,Y_{i3},Y_{i4}}{\mbox{\scriptsize{\eqref{case-2a} holds}}}} \mcal{N} 
&\ll \sum_{\colt{Y_1,Y_{i3},Y_{i4}}{\mbox{\scriptsize{\eqref{case-2a} holds}}}}
Y_1Y_{13}Y_{14}Y_{23}Y_{24}Y_{33}Y_{34}\\
&\ll \sum_{\colt{Y_1,Y_{03},Y_{04},Y_{33}}{Y_{14},Y_{23},Y_{24},Y_{34}}} 
Y_1^{1/2}Y_{03}^{1/4}Y_{04}^{3/4}Y_{14}Y_{23}^{1/2}Y_{24}^{3/2}Y_{33}^{3/4}Y_{34}^{5/4}\\
&\ll B^{1/2}\sum_{\colt{Y_1,Y_{03},Y_{04}}{Y_{23},Y_{24},Y_{33},Y_{34}}} 
Y_{03}^{1/4}Y_{04}^{1/4}Y_{23}^{1/2}Y_{24}^{1/2}Y_{33}^{3/4}Y_{34}^{3/4}.
\end{align*}
Now it follows from \eqref{case2} that $Y_{33}\ll
Y_{03}Y_{04}/Y_{34}$.  We may therefore combine this with the first
inequality in \eqref{r1} to conclude that
\begin{align*}
\sum_{\colt{Y_1,Y_{i3},Y_{i4}}{\mbox{\scriptsize{\eqref{case-2a} holds}}}} \mcal{N} 
&\ll B^{1/2}\sum_{\colt{Y_1,Y_{03},Y_{04}}{Y_{23},Y_{24},Y_{34}}} 
Y_{03}Y_{04}Y_{23}^{1/2}Y_{24}^{1/2}
\ll B (\log B)^5,
\end{align*}
which is also satisfactory for the theorem.

Finally we suppose that \eqref{case3} holds.
On combining \eqref{mendip} with the fact that  $Y_{33}Y_{34} \ll
Y_{03}Y_{04}$, we obtain
$$
Y_{33} \ll \min\left\{  \frac{Y_{04}Y_{14}^2Y_{24}^2Y_{34}}{Y_{03}Y_{13}^2Y_{23}^2},
\frac{Y_{03}Y_{04}}{Y_{34}}\right\} \ll
\frac{Y_{04}Y_{14}Y_{24}}{Y_{13}Y_{23}}.
$$
Summing \eqref{N1} over $Y_{33}$ first, with 
$\min\{Y_{03}Y_{04},Y_{33}Y_{34}\}\leq
Y_{03}^{1/2}Y_{04}^{1/2}Y_{33}^{1/2}Y_{34}^{1/2}$, we therefore obtain
\begin{align*}
\sum_{\colt{Y_1,Y_{i3},Y_{i4}}{\mbox{\scriptsize{\eqref{case3} holds}}}} \mcal{N} 
&\ll 
\sum_{\colt{Y_1,Y_{03},Y_{04},Y_{13}}{Y_{14},Y_{23},Y_{24},Y_{34}}} 
Y_1Y_{03}^{1/2}Y_{04}Y_{13}^{1/2}Y_{14}^{3/2}Y_{23}^{1/2}Y_{24}^{3/2}Y_{34}^{1/2}.
\end{align*}
But then  we may sum over $Y_{03}, Y_{13}$ satisfying the inequalities
in \eqref{r1},  and then $Y_1$ satisfying the second inequality in
\eqref{r2}, in order to conclude that
\begin{align*}
\sum_{\colt{Y_1,Y_{i3},Y_{i4}}{\mbox{\scriptsize{\eqref{case3} holds}}}} \mcal{N} 
&\ll B^{1/4} 
\sum_{\colt{Y_1,Y_{04},Y_{13}}{Y_{14},Y_{23},Y_{24},Y_{34}}} 
Y_1Y_{04}^{1/2}Y_{13}^{1/2}Y_{14}^{3/2}Y_{23}^{1/4}Y_{24}^{5/4}Y_{34}^{1/2}\\
&\ll B^{1/2} 
\sum_{\colt{Y_1,Y_{04},Y_{14}}{Y_{23},Y_{24},Y_{34}}} 
Y_1^{1/2}Y_{04}^{1/2}Y_{14}Y_{24}Y_{34}^{1/2}
\ll B(\log B)^5.
\end{align*}
This too is satisfactory for Theorem \ref{4-2}, and thereby completes its proof.

\section{Open problems}

We close this survey article with a list of five open problems relating to
Manin's conjecture for del Pezzo surfaces. In order to encourage
activity we have deliberately selected an array of very concrete problems.

\begin{enumerate}

\item
{\em Establish \eqref{c2} for a non-singular del Pezzo surface of
  degree $4$.}\\
The surface $x_0x_1-x_2x_3=x_0^2+x_1^2+x_2^2-x_3^2-2x_4^2=0$ has
Picard group of rank $5$.

\item
{\em Establish \eqref{c2} for a non-rational del Pezzo surface.}\\
The surface $x_0x_1-x_2^2=x_0x_2-x_1x_2+x_3^2+x_4^2= 0$, which is
isomorphic (over $\ov{\Q}$) to the type \texttt{iii} surface in the
table, is an example of an {\em Iskovskih surface}.  It is not rational
over $\Q$ \cite[Proposition 7.7]{c-t}.

\item
{\em Break the $4/3$-barrier for a non-singular cubic surface.}\\ 
We have yet to prove an upper bound of the shape
$N_{U,H}(B)=O_S(B^{\theta})$, with $\theta<4/3$, for a single
non-singular cubic surface $S\subset \bfP^3$.   Of course
the ultimate goal is to do this for every such surface, but this
seems to be much harder when the surface doesn't have a conic
bundle structure over $\Q$.  
The surface $x_0x_1(x_0+x_1)=x_2x_3(x_2+x_3)$ admits such a structure
---can one break the $4/3$-barrier for this example?

\item
{\em Establish the lower bound $\NU(B)\gg B (\log B)^3$ for the Fermat
  cubic.}\\ 
The Fermat cubic $x_0^3+x_1^3=x_2^3+x_3^3$ has Picard group of rank $4$.


\item
{\em Better bounds for del Pezzo surfaces of degree $2$}.\\
The arithmetic of non-singular  del Pezzo surfaces of degree $2$ is still very
elusive. These surfaces 
take the shape $t^2=F(x_0,x_1,x_2)$ for a non-singular quartic form
$F$. Let $N(F;B)$ denote the number of integers $t,x_0,x_1,x_2$ such that $t^2=F(\x)$
and $|\x|\leq B$. Can one prove that we always have 
$N(F;B)=O_{\ve,F}(B^{2+\ve})$?  Such an estimate would be essentially
best possible, as consideration of the form $F_0(\x)=x_0^4+x_1^4-x_2^4$ shows.
The best result in this direction is
due to Broberg \cite{bro}, who has established the weaker bound $N(F;B)=O_{\ve,F}(B^{9/4+\ve})$.
For certain quartic forms, such as $F_1(\x)=x_0^4+x_1^4+x_2^4$, 
the Manin conjecture implies that one ought to be able to replace the exponent $2+\ve$ by $1+\ve$. 
Can one prove that $N(F_1;B)=O(B^{\theta})$ for some $\theta<2$?

\end{enumerate}

\begin{ack}
The author is extremely grateful to R\'egis de la Bret\`eche and Per
Salberger, who have both made a number of useful comments about an
earlier version of this paper.
\end{ack}

\end{document}